\documentclass{article}

\usepackage{amsmath,amsthm,amssymb}
\usepackage{graphicx,psfrag,epsfig,multirow}

\begin{document}
\title {A stretched exponential bound on the rate of growth
of the number of periodic points for prevalent diffeomorphisms}
\author{Vadim Yu. Kaloshin \& Brian R. Hunt}
\date{}
\maketitle

\markboth{Periodic points}{Vadim Yu. Kaloshin \& Brian R. Hunt}

\theoremstyle{plain}
\newtheorem{Thm}{Theorem}[section]
\newtheorem{Lm}{Lemma}[section]
\newtheorem{Prop}{Proposition}[section]
\newtheorem{Not}{Notation}[section]
\newtheorem{Cor}{Corollary}[section]
\newtheorem{Def}{Definition}[section]
\newtheorem{Rm}{Remark}[section]
\newtheorem{Q}{Question}
\newtheorem{Prob}{Problem}

\def\cal{\mathcal}
\def\bdef{\begin{Def}}
\def\endef{\end{Def}}
\def\bthm{\begin{Thm}}
\def\ethm{\end{Thm}}
\def\blm{\begin{Lm}}
\def\elm{\end{Lm}}
\def\bprop{\begin{Prop}}
\def\eprop{\end{Prop}}
\def\bcor{\begin{Cor}}
\def\ecor{\end{Cor}}
\def\brm{\begin{Rm}}
\def\erm{\end{Rm}}
\def\bprob{\begin{Prob}}
\def\eprob{\end{Prob}}
\def\beal{\begin{aligned}}
\def\enal{\end{aligned}}
\def\diffeo{dif\-fe\-o\-mor\-phism\ }
\def\beq{\begin{eqnarray}}
\def\eneq{\end{eqnarray}}
\def\Cal{\mathcal}
\def\pseudo{pseudotrajectory\ }
\def\nbhd{ne\-i\-gh\-bor\-ho\-od}
\def\hyp{hy\-per\-bo\-li\-ci\-ty}
\def\phyp{pseudo\-hy\-per\-bo\-li\-ci\-ty}
\def\evg{eigenvalue}
\def\evs{eigenvalues}
\def\pa{\partial}
\def\R{\mathbb R}
\def\Z{\mathbb Z}
\def\L{\Cal L}
\def\J{\Cal J}
\def\bJ{\mathbb J}
\def\D{\Cal D}
\def\P{\mathcal P}
\def\dt{\delta}
\def\al{\alpha}
\def\r{\mathbf r}
\def\gm{\gamma}
\def\gmn{\gm_n}
\def\tgmn{\tilde \gmn}
\def\gmncdt{\gm_n(C,\dt)}
\def\lb{\lambda}
\def\sg{\sigma}
\def\eps{\varepsilon}
\def\inv{^{-1}}
\def\di{dist}
\def\X{\mathbf X}
\def\~{\tilde}
\def\veps{{\vec \eps}}
\def\r{\mathbf r}
\def\pa{\partial}

\section{Introduction}
%\subsection{A problem of the growth of number of periodic points
%and decay of hyperbolicity for generic diffeomorphisms.}

Let $\textup{Diff}^r(M)$ be the space of $C^r$ diffeomorphisms  of
a finite-dimensional smooth compact manifold $M$ with the uniform
$C^r$-topology, where $\dim M \geq 2,$ and let $f \in
\textup{Diff}^r(M)$.  Consider the number of periodic points of
period $n$
\beq \label{grow}
P_n(f)=\# \{x \in M:\ \  x=f^n(x)\}.
\eneq
The main question of this paper is:
\begin{Q}
How quickly can $P_n(f)$ grow with $n$ for a ``generic'' $C^r$
diffeomorphism $f$?
\end{Q}
We put the word ``generic'' in brackets because as the reader
will see the answer depends on notion of genericity.

For technical reasons one sometimes counts only {\em isolated\/}
points of period $n$; let
\beq
P^i_n(f)=\# \{x \in M:\ \  x=f^n(x) \textup{ and } y \neq f^n(y)
\textup{ for } y \neq x \textup{ in some neighborhood of } x\}.
\eneq
We call a diffeomorphism $f \in \textup{Diff}^r(M)$ an
{\em Artin-Mazur diffeomorphism\/} (or simply {\em A-M diffeomorphism\/})
if the number of isolated periodic orbits of $f$ grows   at most
exponentially fast, {\it i.e.\/} for some number $C>0$,
\beq
P^i_n(f) \leq \exp(Cn) \quad {\textup{for all }} n \in \mathbb Z_+.
\eneq
Artin \& Mazur {\cite {AM}} proved the following result.
\bthm
For $0\leq r\leq \infty$, A-M diffeomorphisms are dense
in Diff$^r(M)$ with the uniform $C^r$-topology.
\ethm

We say that a point $x \in M$ of period $n$ for $f$ is hyperbolic if
$df^n(x)$, the derivative of $f^n$ at $x$, has no eigenvalues with
modulus $1$.  (Notice that a hyperbolic solution to $f^n(x) = x$ must
also be isolated.)  We call $f \in \textup{Diff}^r(M)$ a strongly
Artin-Mazur diffeomorphism if for some number $C > 0$,
\beq
P_n(f) \leq \exp(Cn) \quad {\textup{for all }} n \in \mathbb Z_+,
\eneq
and all periodic points of $f$ are hyperbolic (whence $P_n(f) =
P^i_n(f)$).  In {\cite{K1}} an elementary proof of the following
extension of the Artin-Mazur result is given.
\bthm
For $0\leq r<\infty$, strongly A-M diffeomorphisms are dense in
Diff$^r(M)$ with the uniform $C^r$-topology.
\ethm

According to the standard terminology, a set in Diff$^r(M)$ is
called residual if it contains a countable intersection of open
dense sets and a property is called (Baire) generic if
diffeomorphisms with that property form a residual set.
It turns out the A-M property is not generic, as is shown
in {\cite{K2}}. Moreover:
\bthm {\cite{K2}} \label{supergrowth} For any $2\leq r <\infty$
there is an open set $\Cal N \subset$ Diff$^r(M)$ such that for
any given sequence $a=\{a_n\}_{n \in \mathbb Z_+}$ there is a Baire
generic set $\cal R_a$ in $\Cal N$ depending on the sequence $a_n$
with the property if $f \in \cal R_a$, then for infinitely many
$n_k \in \mathbb Z_+$ we have $P^i_{n_k}(f)>a_{n_k}$.
\ethm
\noindent
Of course since $P_n(f) \geq P^i_n(f)$, the same statement can be made
about $P_n(f)$.  But in fact it is shown in \cite{K2} that $P_n(f)$ is
infinite for $n$ sufficiently large, due to a continuum of periodic
points, for at least a dense set of $f \in \cal N$.

The proof of this Theorem is based on a result of
Gonchenko-Shilnikov-Turaev \cite{GST1}. Two slightly
different detailed proofs of their result are given
in {\cite{K2}} and {\cite{GST2}}. The proof in {\cite{K2}}
relies on a strategy outlined in \cite{GST1}.

However, it seems unnatural that if you pick a diffeomorphism
at random then it may have an arbitrarily fast growth of
number of periodic points. Moreover, Baire generic sets in
Euclidean spaces can have zero Lebesgue measure. Phenomena
that are Baire generic, but have a small probability are
well-known in dynamical systems, KAM theory, number theory,
etc.\ (see {\cite{O}}, {\cite{HSY}}, {\cite{K3}} for various
examples). This partially motivates the problem posed by
Arnold {\cite{A}}:
\begin{Prob}
Prove that ``with probability one'' $f \in \textup{Diff}^r(M)$ is an
A-M diffeomorphism.
\end{Prob}

Arnold suggested the following interpretation of ``with probability
one'': {\em for a (Baire) generic finite parameter family of
diffeomorphisms $\{f_\eps\}$, for Lebesgue almost every $\eps$ we have
that $f_\eps$ is A-M\/} ({\it cf.\/} \cite{K3}).  As Theorem 1.3 shows,
a result on the genericity of the set of A-M diffeomorphisms based on
(Baire) topology is likely to be extremely subtle, if possible at
all\footnote{ For example, using techniques from \cite{GST2} and
\cite{K2} one can prove that for a Baire generic finite-parameter
family $\{f_\eps\}$ and a Baire generic parameter value $\eps$
the corresponding diffeomorphism $f_\eps$ is not A-M. Unfortunately,
how to estimate from below the measure of non-A-M diffeomorphisms
in a Baire generic finite-parameter family is
so far an unreachable question.}.  We use instead
a notion of ``probability one'' based on prevalence \cite{HSY,K3},
which is independent of Baire genericity.  We also are able to state
the result in the form Arnold suggested for generic families using this
measure-theoretic notion of genericity.

For a rough understanding of prevalence, consider a Borel measure
$\mu$ on a Banach space $V$.  We say that a property holds
``$\mu$-almost surely for perturbations'' if it holds on a Borel set
$P \subset V$ such that {\it for all $v \in V$ we have $v + w \in P$
for almost every $w$ with respect to $\mu$}\footnote{A similar notion
of prevalence is used in \cite{VK}.}.  Notice that if $V = \mathbb
R^k$ and $\mu$ is Lebesgue measure, then ``$\mu$-almost surely for
perturbations'' is equivalent to ``Lebesgue almost
everywhere''.  Moreover, the Fubini/Tonelli Theorem implies that if
$\mu$ is any Borel probability measure on $\mathbb R^k$, then a
property that holds $\mu$-almost surely for perturbations
must also be hold Lebesgue almost everywhere.  Based on this
observation, we call a property on a Banach space ``prevalent'' if it
holds $\mu$-almost surely for perturbations for some
Borel probability measure $\mu$ on $V$, which for technical reasons
({\it cf.\/} \cite{HSY}) we
require to have compact support.  In order to apply this notion to the
Banach manifold $\textup{Diff}^r(M)$, we must describe how we make
perturbations in this space, which we will do in the next Section.

Our first main result is a partial solution to
Arnold's problem. It says that {\em for a prevalent diffeomorphism
$f\in$ \textup{Diff}$^r(M)$, with $1 < r \leq \infty$, and all $\dt>0$
there exists $C=C(\dt)>0$ such that for all $n \in \Z_+$,\/}
\beq \label{growthbound}
P_n(f) \leq \exp(Cn^{1+\dt}).
\eneq

The Kupka-Smale theorem (see {\it e.g.\/} {\cite{PM}}) states that for
a generic diffeomorphism all periodic points are hyperbolic and
all associated stable and unstable manifolds intersect one another
transversally.  Ref.~{\cite {K3}} shows that the Kupka-Smale theorem
also holds on a prevalent set.
So, the Kupka-Smale theorem, in particular, says that a Baire generic
(resp.\ prevalent) diffeomorphism has only hyperbolic periodic points,
but {\em how hyperbolic are the periodic points, as function of their
period, for a Baire generic (resp.\ prevalent) diffeomorphism $f$?\/}
This is the second main problem we deal with in this paper.

Recall that a linear operator $L:\R^N \to \R^N$ is
{\em hyperbolic\/} if it has no eigenvalues on the unit circle
$\{|z|=1\}\subset \mathbb C$. Denote by $|\cdot|$ the Euclidean
norm in $\mathbb C^N$. Then we define the {\em hyperbolicity\/} of
a linear operator $L$ by
\beq \label{hyp-lin}
\gm(L)=\inf_{\phi \in [0,1)}\inf_{|v|=1}|Lv - \exp(2\pi i \phi) v|.
\eneq
We also say that $L$ is $\gm$-hyperbolic if $\gm(L)\geq\gm$.
In particular, if $L$ is $\gm$-hyperbolic, then its eigenvalues
$\{\lb_j\}_{j=1}^N\subset \mathbb C$ are at least $\gm$-distant
from the unit circle, {\it i.e.\/} $\min_j ||\lb_j|-1|\geq \gm$.
The {\em hyperbolicity\/} of a periodic point $x=f^n(x)$ of period
$n$, denoted by $\gm_n(x,f)$, equals the hyperbolicity of
the derivative $df^n(x)$ of $f^n$ at points $x$, {\it i.e.\/}
$\gm_n(x,f)=\gm(df^n(x))$. Similarly to the number of periodic
points $P_n(f)$ of period $n$, define
\beq \label{hyp-map}
\gm_n(f)=\min_{\{x:\ x=f^n(x)\}}\gm_n(x,f).
\eneq

The idea of Gromov {\cite{G}} and Yomdin {\cite{Y}} of measuring
hyperbolicity is that a $\gm$-hyperbolic point of period $n$ of a
$C^2$ diffeomorphism $f$ has an $M_2^{-2n}\gm$-neighborhood (where
$M_2 = \|f\|_{C^2}$) free from periodic points of the same
period\footnote{In {\cite{Y}} hyperbolicity is introduced as the
minimal distance of eigenvalues to the unit circle. This way of
defining hyperbolicity does not guarantee the existence of a
$M_2^{-2n}\gm$-neighborhood free from periodic points of the same
period \cite{KH}.}. One can prove the following slightly more
general result.

\bprop \label{per-hyp} Let $M$ be a compact manifold of dimension $N$,
let $f:M \to M$ be a $C^{1+\rho}$ diffeomorphism (where $0 < \rho
\leq 1$) that has only hyperbolic periodic points, and let
$M_{1+\rho} = \max(\|f\|_{C^{1+\rho}}, 2^{1/\rho})$.  Then there is
a constant $C=C(M)>0$ such that for each $n\in \Z_+$ we have
\beq
P_n(f)\leq C\ \left(M_{1+\rho}\right)^{nN(1 +
\rho)/\rho}\gm_n(f)^{-N/\rho}.
\eneq
\eprop

Proposition \ref{per-hyp} implies that a lower estimate on
the decay of hyperbolicity $\gm_n(f)$ gives an upper estimate
on the growth of the number of periodic points $P_n(f)$.  Therefore, a natural
question is:
\begin{Q}
How quickly can $\gm_n(f)$ decay with $n$ for a ``generic'' $C^r$
diffeomorphism $f$?
\end{Q}
The existence of lower bound on a rate of decay of $\gm_n(f)$ for
Baire generic $f \in \textup{Diff}^r(M)$ would imply
the existence of an upper bound on a rate of growth of the number
of periodic points $P_n(f)$, whereas no such bound exists by
Theorem \ref{supergrowth}. Thus again we consider genericity in
the measure-theoretic sense of prevalence.  Our second main result,
which in view of Proposition \ref{per-hyp} implies the first main
result, is that {\em for a prevalent diffeomorphism
$f\in$ \textup{Diff}$^r(M)$, with $1 < r \leq \infty$, and all
$\dt>0$ there exists $C=C(\dt)>0$ such that\/}
\beq \label{hyperbolicitydecay}
\gm_n(f) \geq \exp(-Cn^{1+\dt}).
\eneq

Now we shall discuss in more detail our definition of prevalence
(``probability one'') in the space of diffeomorphisms
$\textup{Diff}^r(M)$.

\section{
Prevalence in the Space of Diffeomorphisms
$\textup{Diff}^r(M)$}\label{prevalence}

The space of $C^r$ diffeomorphisms $\textup{Diff}^r(M)$
of a compact manifold $M$ is a Banach manifold.  Locally we can
identify it with a Banach space, which gives it a local linear
structure in the sense that we can perturb a diffeomorphism by
``adding'' small elements of the Banach space.  As we described in
the previous section, the notion of prevalence requires us to make
additive perturbations with respect to a probability measure that
is independent of the place that we make the perturbation.  Thus
although there is not a unique way to put a linear structure on
$\textup{Diff}^r(M)$, it is important to make a choice that is
consistent throughout the Banach manifold.

The way we make perturbations on $\textup{Diff}^r(M)$ by small
elements of a Banach space is as follows.  First we embed $M$ into
the interior of the closed unit ball $B^N \subset \R^N$, which we
can do for $N$ sufficiently large by the Whitney Embedding Theorem
\cite{W}. We emphasize that our results hold for {\it{every}} possible
choice of an embedding of $M$ into $\R^N$.  We then consider a closed
tube neighborhood $U \subset B^N$ of $M$ and the Banach space $C^r(U,\R^N)$
of $C^r$ functions from $U$ to $\R^N$.  Next, we extend every element
$f \in \textup{Diff}^r(M)$ to an element $F \in C^r(U,\R^N)$ that is
strongly contracting in the directions transverse to $M$.  Again the
particular choice of how we make this extension is {\it{not}} important
to our results; in the Appendix  we describe the conditions we
need to ensure that the results of Sacker \cite{Sac} and Fenichel
\cite{F} apply as follows.  Since $F$ has $M$ as an invariant
manifold, if we add to $F$ a small perturbation in $g \in
C^r(U,\R^N)$, the perturbed map $F + g$ has an invariant manifold
in $U$ that is close to $M$.  Then $F + g$ restricted to its
invariant manifold corresponds in a natural way to an element of
$\textup{Diff}^r(M)$, which we consider to be the perturbation of $f
\in \textup{Diff}^r(M)$ by $g \in C^r(U,\R^N)$.  The details of this
construction are described in the Appendix.

In this way we reduce the problem to the study of maps in
$\textup{Diff}^r(U)$, the open subset of $C^r(U,\R^N)$ consisting
of those elements that are diffeomorphisms from $U$ to some subset
of its interior.  The construction we described in the previous
paragraph ensures that the number of periodic points $P_n(f)$ and
their hyperbolicity $\gamma_n(f)$ for elements of $\textup{Diff}^r(M)$
are the same for the corresponding elements of $\textup{Diff}^r(U)$,
so the bounds that we prove on these quantities for almost every
perturbation of any element of $\textup{Diff}^r(U)$ hold as well for
almost every perturbation of any element of $\textup{Diff}^r(M)$.
Another justification for considering diffeomorphisms in Euclidean
space is that the problem of exponential/superexponential growth of
the number of periodic points $P_n(f)$ for a prevalent
$f \in \textup{Diff}^r(M)$ is a {\em local problem\/} on $M$ and
is not affected by a global shape of $M$.

The results stated in the next section apply to any compact domain
$U \subset \R^N$, but for simplicity we state them for the closed unit
ball $B^N$.  In the previous section, we said that a property is
{\em prevalent\/} on a Banach space such as $C^r(B^N,\R^N)$ if
it holds on a Borel subset $S$ for which there exists a Borel
probability measure $\mu$ on $C^r(B^N,\R^N)$ with compact support
such that for all $F \in C^r(B^N,\R^N)$ we have $F + g \in S$ for
almost every $g$ with respect to $\mu$.  The complement of
a prevalent set is said to be {\em shy\/}. We then say that
a property is prevalent on an open subset of $C^r(B^N,\R^N)$ such as
$\textup{Diff}^r(B^N)$ if the exceptions to the property in
$\textup{Diff}^r(B^N)$ form a shy subset of $C^r(B^N,\R^N)$.

In this paper the perturbation measure $\mu$ that we use is supported
within the analytic functions in $C^r(B^N,\R^N)$.  In this sense
we foliate $\textup{Diff}^r(B^N)$ by analytic leaves that are compact
and overlapping.  The main result then says that {\em for every
analytic leaf $L \subset \textup{Diff}^r(B^N)$ and every $\dt>0$,
for almost every diffeomorphism $f \in L$ in the leaf $L$ both
(\ref{growthbound}) and (\ref{hyperbolicitydecay}) are satisfied.\/}
Now we define an analytic leaf as a ``Hilbert brick'' in the space of
analytic functions, and a natural Lebesgue product probability measure
$\mu$ on it.

\section{Formulation of Main Results}\label{mainres}

Fix a coordinate system $x=(x_1,\dots,x_N) \in \R^N \supset B^N$ and
the scalar product $\langle x,y \rangle = \sum_i x_i y_i$.  Let
$\alpha=(\alpha_1,\dots,\alpha_N)$ be a multiindex from $\Z_+^N$, and
let $|\alpha|=\sum_i \alpha_i$.  For a point $x=(x_1,\dots,x_N) \in
\R^N$ we write $x^\alpha=\prod_{i=1}^N x_i^{\alpha_i}$.  Associate to
a real analytic function $\phi:B^N \to \R^N$ the set of coefficients
of its expansion:
\beq \label{expan}
\phi_\veps(x)=\sum_{\alpha \in \Z_+^N}
\veps_\alpha x^\alpha.
\eneq
Denote by $W_{k,N}$ the space of $N$-component homogeneous
vector-polynomials of degree $k$ in $N$ variables and by
$\nu(k,N)=\dim W_{k,N}$ the dimension of $W_{k,N}$. According
to the notation of the expansion (\ref{expan}),
denote coordinates in $W_{k,N}$ by
\beq \label{homogcoord}
\veps_k=
\left(\{\veps_{\alpha}\}_{|\alpha|=k}\right)\in W_{k,N}.
\eneq
In $W_{k,N}$ we use a scalar product that is invariant with respect to
orthogonal transformation of $\R^N \supset B^N$, defined as follows:
\beq \label{scalprod}
\langle \veps_k,\vec \zeta_k \rangle_k=
\sum_{|\al|=k} \binom{k}{\al}^{-1} \langle\veps_{\al},\vec
\zeta_{\al}\rangle, \quad \|\veps_k\|_k=
\bigl(\langle \veps_k,\veps_k \rangle_k \bigr)^{1/2}.
\eneq
Denote by
\beq \label{ball}
B^N_k(r)=\left\{\veps_k \in W_{k,N}: \ \|\veps_k\|_k \leq
r\right\}
\eneq
the closed $r$-ball in $W_{k,N}$ centered at the origin.  Let
$Leb_{k,N}$ be Lebesgue measure on $W_{k,N}$ induced by the scalar
product (\ref{scalprod}) and normalized by a constant so that the
volume of the unit ball is one: $Leb_{k,N}(B^N_k(1))=1$.

Fix a nonincreasing sequence of positive numbers
$\r=\left(\{r_k\}_{k=0}^\infty\right)$ such that
$r_k \to 0$ as $k\to \infty$
and define a Hilbert brick of size $\r$
\beq \label{Brick}
\beal
HB^N(\r) = & \{\veps=\{\veps_{\alpha}\}_{\alpha \in \mathbb Z_+^N}:
\ \textup{for all}\ k \in  \mathbb Z_+, \|\veps_k\|_k\leq r_k\}
\\
= &
B^N_0(r_0) \times B^N_1(r_1) \times \dots \times B^N_k(r_k) \times
\dots \subset
W_{0,N} \times W_{1,N} \times \dots \times W_{k,N} \times \dots .
\enal
\eneq
Define {\em a product Lebesgue probability measure $\mu^N_{\r}$
associated to the Hilbert brick $HB^N(\r)$ of size $\r$\/}
by normalizing for each $k \in \Z_+$ the corresponding Lebesgue
measure $Leb_{k,N}$ on $W_{k,N}$ to the Lebesgue probability
measure on the $r_k$-ball $B^N_k(r_k)$:
\beq \label{measure}
\begin{aligned}
\mu_{k,r}^N= r^{-\nu(k,N)}\ Leb_{k,N}\quad
\textup{and} \quad
\mu^N_{\r}=\times_{k=0}^\infty \mu_{k,r_k}^N.
\end{aligned}
\eneq

\bdef \label{admissible}
Let $f \in \textup{Diff}^r(B^N)$ be a $C^r$
diffeomorphism of $B^N$ into its interior. We call
$HB^N(\r)$ a Hilbert brick of an admissible size
$\r = \left(\{r_k\}_{k=0}^\infty\right)$ with respect to $f$
if:
\begin{itemize}
\item[A)] for each $\veps \in HB^N(\r)$, the corresponding
function $\phi_\veps(x)=\sum_{\al \in \Z_+^N} \veps_\al x^\al$
is analytic on $B^N$;
\item[B)] for each $\veps \in HB^N(\r)$, the corresponding
map $f_\veps(x)=f(x)+\phi_\veps(x)$ is a diffeomorphism
from $B^N$ into its interior, {\it i.e.\/}
$\{f_{\veps}\}_{\veps \in HB^N(\r)} \subset
\textup{Diff}^r(B^N)$;
\item[C)] for all $\dt>0$ and all $C>0$, the sequence
$r_k \exp(C k^{1+\dt})\to \infty$ as $k \to \infty$.
\end{itemize}
\endef

\brm The first and second conditions ensure that the family
$\{f_{\veps}\}_{\veps \in HB^N(\r)}$ lie in an analytic leaf within
the class of diffeomorphisms $\textup{Diff}^r(B^N)$. The third
condition provides us enough freedom to perturb. It is important for
our method to have infinitely many parameters to perturb. If $r_k$'s
decayed too fast to zero it would make our family of
perturbations essentially finite-dimensional.
\erm

An example of an admissible sequence $\r =
(\{r_k\}_{k=0}^\infty)$ is $r_k = \tau/k!$, where $\tau$ depends on
$f$ and is chosen sufficiently small to ensure that condition (B)
holds. Notice that the diameter of $HB^N(\r)$ is then
proportional to $\tau$, so that $\tau$ can be chosen as some
multiple of the distance from $f$ to the boundary of
$\textup{Diff}^r(B^N)$.

\smallskip
\noindent
{\bf Main Theorem.}\ {\it For any $0 < \rho \leq \infty$ and any
$C^{1+\rho}$ diffeomorphism $f \in \textup{Diff}^{1+\rho}(B^N)$,
consider a Hilbert brick $HB^N(\r)$ of an admissible size $\r$
with respect to $f$ and the family of analytic perturbations of $f$
\beq \label{anptb}
\{f_\veps(x)=f(x)+\phi_\veps(x)\}_{\veps \in HB^N(\r)}
\eneq
with the Lebesgue product probability measure $\mu^N_{\r}$
associated to $HB^N(\r)$. Then for every $\dt>0$ and
for $\mu^N_{\r}$-a.e.\ $\veps$ there is $C=C(\veps,\dt)>0$
such that for all $n \in \Z_+$
\beq \label{growth}
\beal
\gm_n(f_\veps) > \ \exp(-Cn^{1+\dt}) \ \ \
P_n(f_\veps) < \ \exp(Cn^{1+\dt}).
\enal
\eneq
}
\brm
The fact that the measure $\mu^N_{\r}$ depends on $f$ does not
conform to our definition of prevalence.  However, we can decompose
$\textup{Diff}^r(B^N)$ into a nested countable union of sets  $\cal
S_j$ that are each a positive distance from the boundary of
$\textup{Diff}^r(B^N)$ and for each $j \in \Z^+$ choose an admissible
sequence ${\r}_j$ that is valid for all $f \in \cal S_j$.  Since
a countable intersection of prevalent subsets of a Banach space is
prevalent \cite{HSY}, the Main Theorem implies the results stated in
terms of prevalence in the introduction.
\erm

In the Appendix we deduce from the Main Theorem the following result.
\bthm \label{genparam} Let
$\{f_\eps\}_{\eps \in B^m}\subset \textup{Diff}^{1+\rho}(M)$
be a generic $m$-parameter family of $C^{1+\rho}$ diffeomorphisms
of a compact manifold $M$ for some $\rho>0$.
Then for every $\dt>0$ and a.e.\ $\eps\in B^m$ there is a constant
$C=C(\eps,\dt)$ such that (\ref{growth}) is satisfied for
every $n\in \Z_+$.
\ethm
\noindent
In the Appendix we also give a precise meaning to the term
{\em generic\/}.

Let us formulate the most general result we shall prove.
\bdef Let $\gm\geq 0$ and
$f \in \textup{Diff}^{1+\rho}(B^N)$ be a $C^{1+\rho}$
diffeomorphism for some $\rho>0$. A point $x \in B^N$
is called $(n,\gm)$-periodic if $\|f^n(x)-x\|\leq \gm$
and $(n,\gm)$-hyperbolic if $\gm_n(x,f)=\gm(df^n(x))\geq \gm$.
\endef
\noindent
(Notice that a point can be $(n,\gm)$-hyperbolic regardless of its
periodicity, but this property is of interest primarily for
$(n,\gm)$-periodic points.)  For positive $C$ and $\dt$ let
$\gmncdt=\exp(-C n^{1+\dt})$.
\bthm \label{main}  Given the hypotheses of the Main Theorem,
for every $\dt>0$ and for $\mu^N_{\r}$-a.e.\ $\veps$ there is
$C=C(\veps,\dt)>0$ such that for all $n \in \Z_+$, every
$(n,\gm^{1/\rho}_n(C,\dt))$-periodic point $x \in B^N$ is
$(n,\gmncdt)$-hyperbolic.  (Here we assume $0 < \rho \leq 1$; in a
space $\textup{Diff}^{1+\rho}(B^N)$ with $\rho > 1$, the statement
holds with $\rho$ replaced by $1$.)
\ethm
\noindent
This result together with Proposition \ref{per-hyp}
implies the Main Theorem,
because every periodic point of period $n$ is $(n,\gm)$-periodic
for all $\gm>0$.
\brm \label{domain}
In the statement of the Main Theorem and Theorem \ref{main}
the unit ball $B^N$ can be replaced by a bounded open set
$U \subset \R^N$.  After scaling, $U$ can be considered as
a subset of the unit ball $B^N$.
\erm

One can define a distance on a compact manifold $M$ and almost
periodic points of diffeomorphisms of $M$. Then one
can cover $M=\cup_i U_i$ by coordinate charts and define
hyperbolicity for almost periodic points using these charts
$\{U_i\}_i$ (see {\cite{Y}} for details). This gives
a precise meaning to the following result.
\bthm Let
$\{f_\eps\}_{\eps \in B^m}\subset \textup{Diff}^{1+\rho}(M)$
be a generic $m$-parameter family of diffeomorphisms
of a compact manifold $M$ for some $\rho>0$.  Then for every $\dt>0$
and almost every $\eps\in B^m$ there is a constant
$C=C(\eps,\dt)$ such that every $(n,\gm^{1/\rho}_n(C,\dt))$-periodic
point $x$ in $B^N$ is $(n,\gmncdt)$-hyperbolic.  (Here again we assume
$0 < \rho \leq 1$, replacing $\rho$ with $1$ in the conclusion if
$\rho > 1$.)
\ethm
The meaning of the term generic is the same as in Theorem
\ref{genparam} and is discussed in the Appendix.

\section{Formulation of the main result in
the $1$-di\-men\-si\-o\-nal case}

The proof of the main result about estimating
the rate of growth of the number of periodic points for
diffeomorphisms in $N$ dimensions has a lot of
complications related to multidimensionality. To describe a model
which is from one side nontrivial and from another side is useful
for understanding the general technique we apply our method for
the $1$-dimensional maps. The statement of the main result for
the $1$-dimensional maps has another important feature: it explains
the statement of the main multidimensional result.

Fix the interval $I=[-1,1]$. Associate to a real analytic function
$\phi:I\to \R$ the set of coefficients of its expansion
\beq
\phi_\eps(x)=\sum_{k=0}^\infty \eps_k x^k.
\eneq
For a nonincreasing sequence of positive numbers
$\r=(\{r_k\}_{k=0}^\infty)$ such that $r_k\to 0$ as
$k \to \infty$ following the multidimensional notations we define
a Hilbert brick of size $\r$
\beq
HB^1(\r)=\{\eps=\{\eps_k\}_{k=0}^\infty:\ \ \
\textup{for all}\ \ \ k\in \Z_+,\ \ \  |\eps_k|\leq r_k\}
\eneq
and the product probability measure $\mu^1_{\r}$ associated to
the Hilbert brick $HB^1(\r)$ of size $\r$ which considers
each $\eps_k$ to be an independent random variable uniformly
distributed on $[-r_k,r_k]$.

\smallskip
\noindent
{\bf Main $1$-dimensional Theorem.}\ {\it For any $0 < \rho \leq \infty$
and any $C^{1+\rho}$ map  $f: I \to I$ of the interval $I=[-1,1]$
consider a Hilbert brick $HB^1(\r)$ of an admissible size
$\r$ with respect to $f$ and the family of analytic perturbations
of $f$
\beq \label{anptb1}
\{f_\eps(x)=f(x)+\phi_\eps(x)\}_{\eps \in HB^1(\r)}
\eneq
with the Lebesgue product probability measure $\mu^1_{\r}$
associated to $HB^1(\r)$. Then for every $\dt>0$ and
$\mu^1_{\r}$-a.e.\ $\eps$ there is $C=C(\eps,\dt)>0$
such that for all $n \in \Z_+$
\beq \label{growth1}
\gm_n(f_\eps) > & \exp(-Cn^{1+\dt}), \quad
P_n(f_\eps) < & \exp(Cn^{1+\dt}).
\eneq
}

In \cite{MMS} Martens-de Melo-Van Strien prove in a sense a stronger
statement for $C^2$ maps. They show that for any $C^2$ map $f$ of
an interval without ``flat'' critical points there are some $\gm>0$ and
$n_0\in \Z_+$ such that for any $n>n_0$ we have $\gm_n(f)>1+\gm$.
This also implies that the number of periodic points is bounded by
an exponential function of the period.  The notion of a flat critical point
used in \cite{MMS} is a nonstandard one from a point of view of
singularity theory.  For $C^2$ maps, they call $x_0$ a flat critical
point of $f$ if $f'(x_0) = f''(x_0) = 0$; the distance from $f(x)$ to
$f(x_0)$ does not have to decay to $0$ as $x\to x_0$ faster than any
power of $x - x_0$.

In \cite{KK} an example of a $C^2$-unimodal map with a critical
point having tangency of order $4$ and an arbitrarily fast rate of
growth of the number of periodic points is given.
Another advantage of the Main $1$-dimensional Theorem is that it works
for $C^{1+\rho}$ maps with $0< \rho<1$, whereas the result
in \cite{MMS} works only for $C^2$ maps.

\section{
Strategy of the proof}\label{strategy}
Here we describe the strategy of the proof of Theorem \ref{main}.  The
basic technique is developed and many of the technical difficulties
are resolved by the first author in \cite{K4}.  The
general idea is to fix $C > 0$ and prove an upper bound on the
measure of the set of ``bad'' parameter values $\veps \in HB^N(\r)$
for which the conclusion of the theorem does not hold.  The upper
bound we obtain will approach zero as $C \to \infty$, from which it
follows immediately that the set of $\veps \in HB^N(\r)$ that are
``bad'' for all $C > 0$ has measure zero.  For a given $C > 0$, we
bound the measure of ``bad'' parameter values inductively as follows.

{\it Stage 1}.\ We delete all parameter values $\veps\in
HB^N(\r)$ for which the corresponding diffeomorphism $f_\veps$
has an almost fixed point that is not sufficiently
hyperbolic, and bound the measure of the deleted set.

{\it Stage 2}.\ After Stage 1, we consider only parameter values for which
all almost fixed points are sufficiently hyperbolic. Then we delete
all parameter values $\veps$ \ for which
$f_\veps$ has an almost periodic point of period $2$ which is not
sufficiently hyperbolic and bound the measure of that set.

{\it Stage n}.\ We consider only parameter values for which
all almost periodic point of period at most $n-1$ are sufficiently
hyperbolic (we shall call this {\em the Inductive Hypothesis\/}). Then
we delete all parameter values $\veps$ for which
$f_\veps$ has an almost periodic point of period $n$ which
is not sufficiently hyperbolic and bound the measure of that set.

The main difficulty in the proof is then to find a bound on the
measure of ``bad'' parameter values at stage $n$ such that
the bounds are summable over $n$ and that the sum approaches zero as
$C \to \infty$.  Let us formalize the problem.  Fix positive $\rho$,
$\dt$, and $C$, and recall that $\{\gmncdt=\exp(-C n^{1+\dt})\}$ for
$n \in \Z_+$.  Assume $\rho \leq 1$; if not, change its value to $1$.
\bdef \label{inductive}
A diffeomorphism $f\in \textup{Diff}^{1+\rho}(B^N)$
satisfies the Inductive Hypothesis of order $n$ with
constants $(C,\dt,\rho)$, denoted $f\in IH(n,C,\dt,\rho)$, if for all
$k\leq n$, every $(k,\gm^{1/\rho}_k(C,\dt))$-periodic point is
$(k,\gm_k(C,\dt))$-hyperbolic.
\endef
\noindent
For $f \in \textup{Diff}^{1+\rho}(M)$, consider the sequence of sets
in the parameter space
$HB^N(\r)$
\beq \label{badset}
B_n(C,\dt,\rho,\r,f) = \{\veps \in HB^N(\r):
f_\veps \in IH(n-1,C,\dt,\rho)
\textup{ but } f_\veps \notin IH(n,C,\dt,\rho) \}
\eneq
In other words, $B_n(C,\dt,\rho,\r,f)$
is the set of ``bad'' parameter values
$\veps \in HB^N(\r)$ for which all almost periodic
points of $f_\veps$ with period strictly less than $n$ are
sufficiently hyperbolic, but there is an almost periodic point of
period $n$ that is not sufficiently hyperbolic.  Let
\beq \label{norms}
M_1=\sup_{\veps \in HB^N(\r)} \max\{\|f_\veps\|_{C^1},
\|f_\veps\inv\|_{C^1}\};\
M_{1+\rho}=\sup_{\veps \in HB^N(\r)}
\max\{\|f_\veps\|_{C^{1+\rho}},M_1,2^{1/\rho}\}.
\eneq

Our goal is to find an upper bound $\mu_n(C,\dt,\rho,
\r,M_{1+\rho})$ for the measure $\mu_{\r}^N\left(
B_n(C,\dt,\rho,\r,f)\right)$ of the set of ``bad'' parameter
values.  Then $\sum_{n=1}^\infty \mu_n(C,\dt,\rho,\r,M_{1+\rho})$
is an upper bound on the measure of
$\cup_{n=1}^\infty B_n(C,\dt,\rho,\r,f)$, which is the set of
all parameter values $\veps$ for which $f_\veps$ has for some $n$ an
$(n,\gm^{1/\rho}_n(C,\dt))$-periodic point that is not
$(n,\gm_n(C,\dt))$-hyperbolic.  If this sum converges and
\beq \label{converge}
\sum_{n=1}^\infty \mu_n(C,\dt,\rho,\r,M_{1+\rho})=
\mu(C,\dt,\rho,\r,M_{1+\rho})\to 0 \ \textup{as}\ C \to \infty
\eneq
for every positive $\rho$, $\dt$, and $M_{1+\rho}$, then Theorem
\ref{main} follows.  In the remainder of this announcement we describe
the key construction we use to obtain a bound $\mu_n(C,\dt,\rho,
\r,M_{1+\rho})$ that meets condition (\ref{converge}).

\section{
Perturbation of recurrent trajectories by
Lagrange interpolation polynomials}\label{Lagrangestart}

The approach we take to estimate the measure of ``bad'' parameter
values in the space of perturbations $HB^N(\r)$ is to choose a
coordinate system for this space and for a finite subset of the
coordinates to estimate the amount that we must change a particular
coordinate to make a ``bad'' parameter value ``good''.  Actually we
will choose a coordinate system that depends on a particular point
$x_0 \in B^N$, the idea being to use this coordinate system to
estimate the measure of ``bad'' parameter values corresponding to
initial conditions in some neighborhood of $x_0$, then cover $B^N$
with a finite number of such neighborhoods and sum the corresponding
estimates.  For a particular set of initial conditions, a
diffeomorphism will be ``good'' if every point in the set is either
sufficiently nonperiodic or sufficiently hyperbolic.

In order to keep the notations and formulas simple as we formalize
this approach, we consider the case of 1-dimensional maps, but the
reader should always have in mind that our approach is designed for
multidimensional diffeomorphisms.  Let $f:I \to I$ be a
$C^1$ map on the interval $I=[-1,1]$. Recall that a trajectory
$\{x_k\}_{k \in \Z}$ of $f$ is called {\em recurrent\/} if
it returns arbitrarily close to its initial position ---
that is, for all $\gm>0$ we have $|x_0-x_n|<\gm$ for some $n>0$.
A very basic question is how much one should perturb $f$ to
make $x_0$ periodic.  Here is an elementary Closing Lemma that
gives a simple partial answer to this question.

\smallskip
\noindent
{\bf Closing Lemma.}\ {\it Let $\{x_k=f^k(x_0)\}_{k=0}^{n}$ be
a trajectory of length $n+1$ of a map $f:I \to I$.
Let $u=(x_0-x_n)/\prod_{k=0}^{n-2}(x_{n-1}-x_k)$.
Then $x_0$ is a periodic point of period $n$ of the map
\beq \label{closing}
f_u(x)=f(x)+u\prod_{k=0}^{n-2}(x-x_k).
\eneq}
\noindent
Of course $f_u$ is close to $f$ if and only if $u$ is
sufficiently small, meaning that $|x_0-x_n|$ should be small
compared to $\prod_{k=0}^{n-2}|x_{n-1}-x_k|$. However,
this product is likely to contain small factors for
recurrent trajectories. In general, it is difficult to control
the effect of perturbations for recurrent trajectories. The simple
reason why is because {\em one can not perturb $f$ at two
nearby points independently\/}.

The Closing Lemma above also gives an idea of how much we must change
the parameter $u$ to make a point $x_0$ that is $(n,\gm)$-periodic not
be $(n,\gm)$-periodic for a given $\gm > 0$, which as we described
above is one way to make a map that is ``bad'' for the initial
condition $x_0$ become ``good''.  To make use of our other alternative
we must determine how much we need to perturb a map $f$ to make a
given $x_0$ be $(n,\gm)$-hyperbolic for some $\gm > 0$.

\smallskip
\noindent
{\bf Perturbation of hyperbolicity.}\ {\it
Let $\{x_k=f^k(x_0)\}_{k=0}^{n-1}$ be a trajectory
of length $n$ of a $C^1$ map $f:I \to I$.  Then for the map
\beq \label{movehyperbolicity}
f_v(x)=f(x)+v(x-x_{n-1})\prod_{k=0}^{n-2}(x-x_k)^2
\eneq
such that $v\in \R$ and
\beq
\left|\vphantom{{f'}^2}|(f^n_v)'(x_0)| - 1\right| =
\left|\vphantom{{{\prod_0^n}^2}^2}
\left|\prod_{k=0}^{n-1}f'(x_k)+
v\prod_{k=0}^{n-2}(x_{n-1}-x_k)^2\prod_{k=0}^{n-2}f'(x_k)\right|
- 1\right| > \gm
\eneq
we have that $x_0$ is an $(n,\gm)$-hyperbolic point of $f_v$.}
\noindent

One more time we can see the product of distances
$\prod_{k=0}^{n-2}|x_{n-1}-x_k|$ along the trajectory
is important quantitative characteristic of how much freedom
we have to perturb.

The perturbations (\ref{closing}) and (\ref{movehyperbolicity}) are
reminiscent of Lagrange interpolation polynomials.  Let us put these
formulas into a general setting using singularity theory.

Given $n > 0$ and a $C^1$ function $f : I \to \R$ we define
an associated function $j^{1,n}f: I^n \to I^n \times \R^{2n}$ by
\beq
\label{multijet}
j^{1,n}f(x_0, \dots, x_{n-1}) = \left(x_0, \dots, x_{n-1}, f(x_0),
\dots, f(x_{n-1}), f'(x_0), \dots, f'(x_{n-1})\right).
\eneq
In singularity theory this function is called the {\em $n$-tuple
$1$-jet\/} of $f$. The ordinary $1$-jet of $f$, usually
denoted by $j^1 f(x) = (x,f(x),f'(x))$, maps $I$ to
the {\em $1$-jet space\/} $\J^1(I,\R) \simeq I \times \R^2$.
The product of $n$ copies of $\J^1(I,\R)$, called the {\em multijet
space\/}, is denoted by
\beq
\J^{1,n}(I,\R)=\underbrace{\J^1(I,\R) \times \dots
\times \J^1(I,\R)}_{n \ \textup{times}},
\eneq
and is equivalent to $I^n \times \R^{2n}$ after rearranging
coordinates.
The $n$-tuple 1-jet of $f$ associates with each $n$-tuple of points
in $I^n$ all the information necessary to determine how close the
$n$-tuple is to being a periodic orbit, and if so, how close it is to
being nonhyperbolic.

The set
\beq
\Delta_n(I)=\left\{\{x_0, \dots, x_{n-1}\}\times \R^{2n} \subset
\J^{1,n}(I,\R) :
\exists\ i \neq j\ \textup{such that}\ x_i=x_j\right\}
\eneq
is called the {\em diagonal\/} (or sometimes the {\em generalized
diagonal\/}) in the space of multijets.
In singularity theory the space of multijets is defined outside
of the diagonal $\Delta_n(I)$ and is usually denoted by
$\J^1_n(I,\R)=\J^{1,n}(I,\R) \setminus \Delta_n(I)$ (see \cite{GG}).
It is easy to see that {\em a recurrent trajectory $\{x_k\}_{k\in \Z_+}$
is located in a neighborhood of the diagonal  $\Delta_n(I)$ in
the space of multijets for a sufficiently large $n$\/}. If
$\{x_k\}_{k=0}^{n-1}$ is a part of a recurrent
trajectory of length $n$, then the product of distances along the
trajectory
\beq \label{productformula}
\prod_{k=0}^{n-2} \left| x_{n-1}-x_k \right|
\eneq
measures how close $\{x_k\}_{k=0}^{n-1}$ is to the diagonal
$\Delta_n(I)$, or how independently one can perturb points of a
trajectory. One can also say that (\ref{productformula})
is a quantitative characteristic of how recurrent a trajectory
of length $n$ is. Introduction of this {\em product of distances
along a trajectory\/} into the analysis of recurrent trajectories
is a new point of our paper.

\section{
Lagrange interpolation and blow-up
along the diagonal in multijet space}\label{blowup}

Now we present a construction due to Grigoriev and Yakovenko
{\cite{GY}} which puts the ``Closing Lemma'' and ``Perturbation of
Hyperbolicity'' statements above into a general framework.
It is an interpretation of Lagrange interpolation
polynomials as an algebraic blow-up along
the diagonal in the multijet space.
In order to keep the notations and formulas simple we continue in this
section to consider only the 1-dimensional case.

Consider the $2n$-parameter family of perturbations of a $C^1$ map
$f:I \to I$ by polynomials of degree $2n-1$
\beq \label{2ndegree}
f_\eps(x)=f(x)+\phi_\eps(x), \qquad \phi_\eps(x) =
\sum_{k=0}^{2n-1}\eps_k x^k,
\eneq
where $\eps = (\eps_0, \ldots, \eps_{2n-1}) \in \R^{2n}$.  The
perturbation vector $\eps$ consists of coordinates from the Hilbert
brick $HB^1(\r)$ of analytic perturbations defined in
Section~\ref{mainres}.
Our goal now is to describe how such perturbations affect the
$n$-tuple $1$-jet of $f$, and since the operator $j^{1,n}$ is linear
in $f$, for the time being we consider only the perturbations
$\phi_\eps$ and their $n$-tuple $1$-jets.  For each $n$-tuple
$\{x_k\}_{k=0}^{n-1}$ there is a natural transformation
$\J^{1,n} : I^n \times \R^{2n} \to \J^{1,n}(I,\R)$ from
$\eps$-coordinates to jet-coordinates, given by
\beq
\label{epstojet}
\J^{1,n}(x_0, \dots, x_{n-1}, \eps) =
j^{1,n} \phi_\eps(x_0, \dots, x_{n-1}).
\eneq

Instead of working directly with the transformation $\J^{1,n}$, we
introduce intermediate $u$-coordinates based on Lagrange interpolation
polynomials.  The relation between $\eps$-coordinates and
$u$-coordinates is given implicitly by
\beq\label{identity}
\phi_\eps(x) = \sum_{k=0}^{2n-1}\eps_k x^k=
\sum_{k=0}^{2n-1}u_k \prod_{j=0}^{k-1} (x-x_{j(\textup{mod}\ n)}).
\eneq
Based on this identity, we will define functions
$\D^{1,n} : I^n \times \R^{2n} \to I^n \times \R^{2n}$
and
$\pi^{1,n}:I^n \times \R^{2n}\to \J^{1,n}(I,\R)$
so that $\J^{1,n} = \pi^{1,n} \circ \D^{1,n}$, or in other words the
diagram in Figure~1 commutes.  We will show later that
$\D^{1,n}$ is invertible, while $\pi^{1,n}$ is invertible away from
the diagonal $\Delta_n(I)$ and defines a blow-up along it in the space
of multijets $\J^{1,n}(I,\R)$.

\begin{figure}[htbp]
  \begin{center}
   \begin{psfrags}
     \psfrag{DivD}{\small{${\cal DD}^{1,n}(I,\R)=
     \underbrace{I \times \dots \times I}_{n \ \textup{times}}
     \times \R^{2n}$}}
     \psfrag{MultJ}{\small{${\cal J}^{1,n}(I,\R)=\underbrace{I \times
     \dots \times I}_{n \ \textup{times}} \times \R^{2n}$}}
     \psfrag{N}{\small{${\D}^{1,n}$}}
     \psfrag{(jf,...,jf)}{\small{${\J}^{1,n}$}}
     \psfrag{Phase}{\small{$\underbrace{I \times \dots
     \times I}_{n \ \textup{times}} \times \R^{2n}$}}
     \psfrag{pi}{\small{$\pi^{1,n}$}}
    \includegraphics[width= 4in,angle=0]{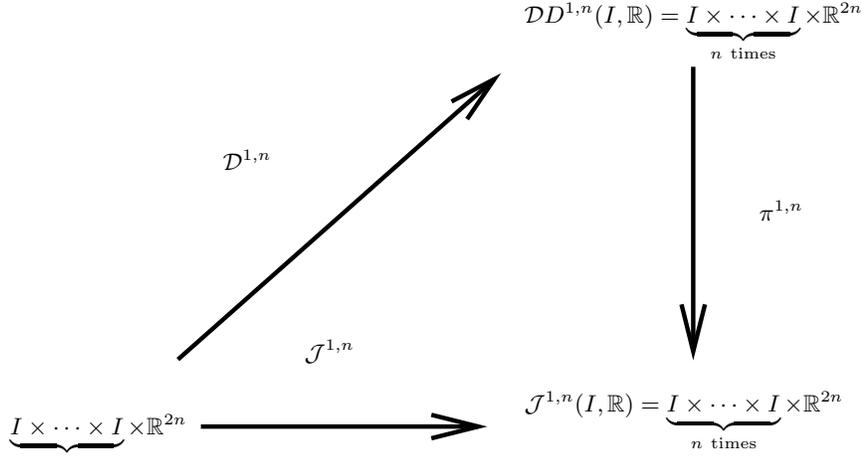}
    \end{psfrags}
   \caption{an Algebraic Blow-up along the Diagonal $\Delta_n(I)$}
  \end{center}
\end{figure}

The intermediate space, which we denote by ${\cal {DD}}^{1,n}(I,\R)$,
is called {\em the space of divided differences\/} and consists of
$n$-tuples of points $\{x_k\}_{k=0}^{n-1}$ and $2n$ real coefficients
$\{u_k\}_{k=0}^{2n-1}$.  Here are explicit coordinate-by-coordinate
formulas defining $\pi^{1,n} : {\cal {DD}}^{1,n}(I,\R) \to
\J^{1,n}(I,\R)$.
\beq
\beal\label{Lagrangeexpression}
\phi_\eps(x_0) = &\, u_0, \\
\phi_\eps(x_1) = &\, u_0+u_1(x_1-x_0),\\
\phi_\eps(x_2) = &\, u_0+u_1(x_2-x_0)+u_2(x_2-x_0)(x_2-x_1),\\
\vdots\,&\\
\phi_\eps(x_{n-1}) = &\, u_0+u_1(x_{n-1}-x_0)+\dots
+u_{n-1}(x_{n-1}-x_0)\dots (x_{n-1}-x_{n-2}),\\
\phi_\eps'(x_0) = &\, \frac{\pa}{\pa x}\left(\sum_{k=0}^{2n-1}
u_k \prod_{j=0}^{k-1} (x-x_{j(\textup{mod}\ n)})\right)\Big|_{x=x_0}, \\
\vdots\,&\\
\phi_\eps'(x_{n-1}) = &\, \frac{\pa}{\pa x} \left(\sum_{k=0}^{2n-1}
u_k \prod_{j=0}^{k-1} (x-x_{j(\textup{mod}\ n)})\right)\Big|_{x=x_{n-1}},
\enal
\eneq

These formulas are very useful for dynamics.  For a given base map $f$
and initial point $x_0$, the image $f_\eps(x_0) = f(x_0) +
\phi_\eps(x_0)$ of $x_0$ depends only on $u_0$.  Furthermore the
image can be set to any desired point by choosing $u_0$ appropriately
--- we say then that it depends nontrivially on $u_0$.  If
$x_0$, $x_1$, and $u_0$ are fixed, the image $f_\eps(x_1)$ of $x_1$
depends only on $u_1$, and as long as $x_0 \neq x_1$ it depends
nontrivially on $u_1$.  More generally for $0 \leq k \leq n-1$, if
pairwise distinct points $\{x_j\}_{j=0}^k$ and coefficients
$\{u_j\}_{j=0}^{k-1}$ are fixed, then the image $f_\eps(x_k)$ of
$x_k$ depends only and nontrivially on $u_k$.

Suppose now that an $n$-tuple of points $\{x_j\}_{j=0}^{n}$ not on the
diagonal $\Delta_n(I)$ and Lagrange coefficients $\{u_j\}_{j=0}^{n-1}$
are fixed.  Then derivative $f'_\eps(x_0)$ at $x_0$ depends only and
nontrivially on $u_n$.  Likewise for $0 \leq k \leq n-1$, if distinct
points $\{x_j\}_{j=0}^{n}$ and Lagrange coefficients
$\{u_j\}_{j=0}^{n+k-1}$ are fixed, then the derivative $f'_\eps(x_k)$
at $x_k$ depends only and nontrivially on $u_{n+k}$.

As Figure~2 illustrates, these considerations show that for any map
$f$ and any desired trajectory of distinct points with any given
derivatives along it, one can choose Lagrange coefficients
$\{u_k\}_{k=0}^{2n-1}$ and explicitly construct a map $f_\eps = f +
\phi_\eps$ with such a trajectory.  Thus we have shown that $\pi^{1,n}$ is
invertible away from the diagonal $\Delta_n(I)$ and defines a blow-up
along it in the space of multijets $\J^{1,n}(I,\R)$.

\begin{figure}[htbp]\label{NIP}
  \begin{center}
   \begin{psfrags}
     \psfrag{x_0}{$x_0$}
     \psfrag{x_1}{$x_1$}
     \psfrag{x_k}{$x_k$}
     \psfrag{x_{k+1}}{$x_{k+1}$}
     \psfrag{f_u(x_0)}{$f_u(x_0)$}
     \psfrag{f_u(x_k)}{$f_u(x_k)$}
     \psfrag{f'_u(x_0)}{$f'_u(x_0)$}
     \psfrag{f'_u(x_k)}{$f'_u(x_k)$}
     \psfrag{u_0}{$u_0$}
     \psfrag{u_k}{$u_k$}
     \psfrag{u_n}{$u_n$}
     \psfrag{u_{n+k}}{$u_{n+k}$}
     \psfrag{dots}{$ \cdots $}
    \includegraphics[width= 4in,angle=0]{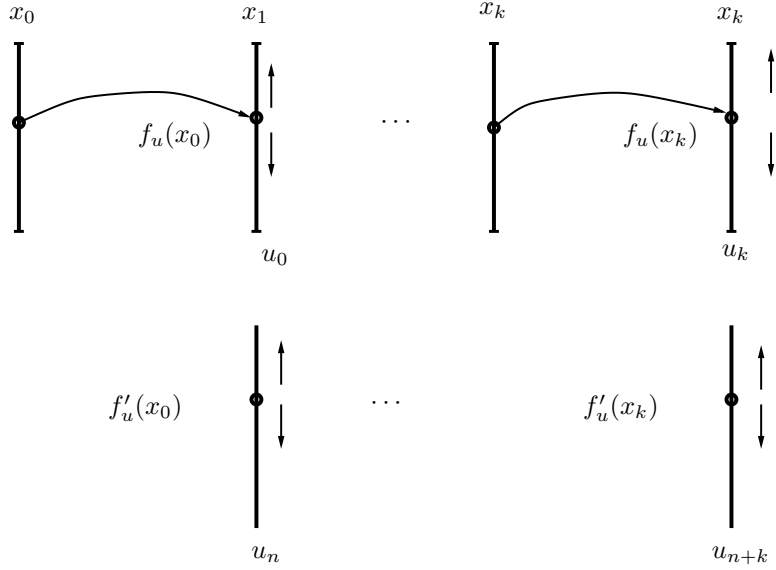}
    \end{psfrags}
   \caption{Lagrange coefficients and their action}
  \end{center}
\end{figure}

Next we define the function $\D^{1,n} : I^n \times \R^{2n} \to
{\cal {DD}}^{1,n}(I,\R)$ explicitly using so-called divided
differences.  Let $g:\R \to \R$ be a $C^r$ function of one real
variable.
\bdef \label{divdif}
The {\em first order divided difference\/} of $g$ is defined as
\beq \begin{aligned}
\Delta g(x_0, x_{1})=\frac {g(x_1)-g(x_0)}{x_1-x_0}
\end{aligned}
\eneq
for $x_1 \neq x_0$ and extended by its limit value as
$g'(x_0)$ for $x_1=x_0$.  Iterating this construction we define
divided differences of the $m$-th order for $2 \leq m \leq r$,
\beq \begin{aligned}
\Delta^m g(x_0, \dots, x_m) =
\frac {\Delta^{m-1} g(x_0, \dots, x_{m-2}, x_m)-
\Delta^{m-1} g(x_0, \dots, x_{m-2}, x_{m-1})}{x_{m}-x_{m-1}}
\end{aligned}
\eneq
for $x_{m-1} \neq x_m$ and extended by its limit value for
$x_{m-1}=x_m$.
\endef
A function loses at most one derivative of smoothness with each
application of $\Delta$, so $\Delta^m g$ is at least $C^{r-m}$ if $g$
is $C^r$.  Notice that $\Delta^m$ is linear as a function of $g$, and
one can show that it is a symmetric function of $x_0,\ldots,x_m$; in
fact, by induction it follows that
\beq
\Delta^m g(x_0, \dots, x_m)=
\sum_{i=0}^m \frac{g(x_i)}{\prod_{j \neq i} (x_i - x_j)}
\eneq
Another identity that is proved by induction will be more important
for us, namely
\beq
\Delta^m\ x^k(x_0,\dots, x_m)=p_{k,m}(x_0,\dots, x_m),
\eneq
where $p_{k,m}(x_0, \dots, x_m)$ is $0$ for $m > k$ and for $m \leq k$
is the sum of all degree $k-m$ monomials in $x_0, \dots, x_m$ with
unit coefficients,
\beq \label{homogpolyn}
p_{k,m}(x_0,\dots, x_m)=\sum_{r_0+\dots + r_m=k-m}
\quad \prod_{j=0}^{m} x_j^{r_j}.
\eneq

The divided differences form coefficients for the Lagrange interpolation
formula.  For all $C^\infty$ functions $g : \R \to \R$ we have
\beq
\beal\label{int}
g(x) = & \Delta^0 g(x_0) +\Delta^1 g(x_0,x_1) (x-x_0)+ \dots \\
& + \Delta^{n-1} g(x_0,\dots, x_{n-1}) (x-x_0) \dots (x-x_{n-2}) \\
& + \Delta^{n} g(x_0,\dots, x_{n-1},x) (x-x_0) \dots (x-x_{n-1})
\enal
\eneq
identically for all values of $x, x_0, \dots, x_{n-1}$.
All terms of this representation are polynomial in $x$
except for the last one which we view as a remainder term.
The sum of the polynomial terms is the degree $(n-1)$ {\em Lagrange
interpolation polynomial\/} for $g$ at $\{x_k\}_{k=0}^{n-1}$.  To
obtain a degree $2n-1$ interpolation polynomial for $g$ and its
derivative at $\{x_k\}_{k=0}^{n-1}$, we simply use (\ref{int}) with
$n$ replaced by $2n$ and the $2n$-tuple of points
$\{x_{k(\textup{mod}\ n)}\}_{k=0}^{2n-1}$.

Recall that $\D^{1,n}$ was defined implicitly by (\ref{identity}).
We have described how to use divided differences to construct a degree
$2n-1$ interpolating polynomial of the form on the right-hand side of
(\ref{identity}) for an arbitrary $C^\infty$ function $g$.  Our
interest then is in the case $g = \phi_\eps$, which as a degree
$2n-1$ polynomial itself will have no remainder term and coincide
exactly with the interpolating polynomial.  Thus $\D^{1,n}$ is given
coordinate-by-coordinate by
\beq
\beal\label{Lagrangemap}
u_m = &\, \Delta^m \left( \sum_{k=0}^{2n-1}\eps_k x^k\right)
(x_0,\dots,x_{m\ (mod\ n)}) \\
= &\, \eps_m+\sum_{k=m+1}^{2n-1} \eps_k p_{k,m}(x_0,\dots , x_{m\ (mod\ n)})
\enal
\eneq
for $m = 0, \dots, 2n-1$.  We call the transformation given by
(\ref{Lagrangemap}) the {\em Lagrange map\/}.  Notice that for fixed
$\{x_k\}_{k=0}^{2n-1}$, the Lagrange map is linear and given by an upper
triangular matrix with units on the diagonal.  Hence it is Lebesgue
volume-preserving and invertible, whether or not
$\{x_k\}_{k=0}^{2n-1}$ lies on the diagonal $\Delta_n(I)$.

We call the basis of monomials
\beq \label{Lagrangebasis}
\prod_{j=0}^k (x-x_{j(\textup{mod}\ n)}) \ \ \
\textup{for}\ \ \ k=0,\dots, 2n-1
\eneq
in the space of polynomials of degree $2n-1$ the {\em Lagrange basis\/}
defined by the $n$-tuple $\{x_k\}_{k=0}^{n-1}$.  The Lagrange map and
the Lagrange basis, and their analogues in dimension $N$, are useful
tools for perturbing trajectories and estimating the measure
$\mu_n(C,\dt,\rho, M_{1+\rho})$
of ``bad'' parameter values $\veps \in HB^N(\r)$.

\section{Discretization method}

The fundamental problem with using the Lagrange basis to estimate the
measure of ``bad'' parameter values, those for which there is an
almost periodic point of period $n$ that is not sufficiently
hyperbolic, is that the Lagrange basis depends on the almost periodic
$n$-tuple $\{x_k\}_{k=0}^{n-1}$.  For a particular ``bad'' parameter
value we can fix this $n$-tuple and the corresponding Lagrange basis,
then estimate the measure of the set of parameters for which a nearby
$n$-tuple is both almost periodic and not sufficiently hyperbolic.
But there are a continuum of possible $n$-tuples, so how can we
account for all of the possible cells of ``bad'' parameter values
$\veps$ within our parameter brick $HB^N(\r)$?  At the beginning
of Section~\ref{Lagrangestart}, we indicated that for a particular
initial condition
$x_0$ we would obtain an estimate on the measure of ``bad'' parameter
values corresponding to an almost periodic point in a neighborhood of
$x_0$, and thus need only to consider a discrete set of initial
conditions.  But as the parameter vector $\veps$ varies over
$HB^N(\r)$, there is (for large $n$ at least) a wide range of
possible length-$n$ trajectories starting from a particular $x_0$, so
there is no hope of using a single Lagrange basis to estimate even the
measure of ``bad'' parameter values corresponding to a single $x_0$.

The solution to this problem is to discretize the entire space of
$n$-tuples $\{x_k\}_{k=0}^{n-1}$, considering only those that lie on a
particular grid.  If we choose the grid spacing small enough, then
every almost periodic orbit of period $n$ that is not sufficiently
hyperbolic will have a corresponding pseudotrajectory of length $n$ on
the grid that also has small hyperbolicity.  In this way we reduce the
problem to bounding the measure of a set of ``bad'' parameter values
corresponding to a particular length $n$ pseudotrajectory, and then
summing the bounds over all possible length $n$ pseudotrajectories on
the chosen grid.

Returning to the general case of $C^{1+\rho}$ diffeomorphisms on
$B^N$, where we assume $0 < \rho \leq 1$, the grid spacing we use at
stage $n$ is $\tgmn(C,\dt,\rho) = N^{-1} (M_{1+\rho}^{-2n}
\gmncdt)^{1/\rho}$, where $M_{1+\rho} > 1$ is a bound on the $C^{1+\rho}$
norm of the diffeomorphisms $f_\veps$ corresponding to parameters
$\veps \in HB^N(\r)$.  This ensures that when rounded off to the
nearest grid points $\{x_k\}_{k=0}^{n-1}$, an almost periodic orbit of
length $n$ becomes an $N(M_{1+\rho} +1)\tgmn(C,\dt,\rho)$-pseudotrajectory,
meaning that $|f_\veps(x_j) -x_{j+1}| \leq N(M_{1+\rho} + 1)
\tgmn(C,\dt,\rho)$ for $j = 0, 1, \ldots,
n-2$.  It also ensures that when rounding, the derivative $df_\veps$
changes by at most $M_{1+\rho}^{1-2n} \gmncdt$, which in turn implies
that the change in hyperbolicity over all $n$ points is small compared
with $\gmncdt$.  (Recall that $\gmncdt$ is our tolerance for
hyperbolicity at stage $n$.)

Roughly speaking, in the case $N = 1$ our estimate on the measure
of ``bad'' parameter values for a particular $n$-tuple
$\{x_k\}_{k=0}^{n-1}$ is then proportional to
$\left(\tgmn(C,\dt,\rho)\right)^n \gmncdt$,
whereas the number of possible $n$-tuples is proportional to
$\left(\tgmn(C,\dt,\rho)\right)^{-n}$, making our bound
$\mu_n(C,\dt,\rho,\r,M_{1+\rho})$
on the total measure of ``bad'' parameter values at stage $N$
proportional to $\gmncdt$.  The remaining problem then is to show that
for maps satisfying the Inductive Hypothesis of order $n-1$, we can
bound the proportionality factor in such a way that
$\mu_n(C,\dt,\rho,\r,M_{1+\rho})$ meets the conditions prescribed
in Section~\ref{strategy}, namely that it be summable over $n$ and
that the sum
approaches $0$ as $C \to \infty$.  (Notice that the sequence $\gmncdt$
meets these conditions.)  The proportionality factor depends on the
product of distances described in Section~\ref{Lagrangestart}, and in
\cite{K4} we proceed as follows.  At the $n$th stage we split length
$n$ trajectories of diffeomorphisms satisfying the Inductive
Hypothesis into three groups.  One group consists of what we call
``simple'' trajectories for which the product of distances is not
too small.  For nonsimple trajectories we show that either
the trajectory is sufficiently hyperbolic by the Inductive Hypothesis
(second group) or the trajectory returns very close to itself before
the $n$th iteration and is simple (not recurrent) up to this point
(third group).  In the latter case, perturbation
by Lagrange polynomials of order lower than $n$ at the point of
a close return can control the behavior of that trajectory up to
length $n$.

Notice that in the preceding paragraph, even if the product of
distances is not small, the proportionality factor in our estimate on
the measure of ``bad'' parameter values for a given $n$-tuple
$\{x_k\}_{k=0}^{n-1}$ is large because the parameter
measure is normalized to be $1$ on a brick $HB^1(\r)$ whose sides
decay rapidly; the normalization increases the measure by a factor of
$r_0 r_1 \cdots r_{n-1} r_{2n-1}$.  However, we are able to
show that when considering only diffeomorphisms $f_\veps$ with $\veps
\in HB^1(\r)$, the number of $n$-tuples we must consider as possible
pseudotrajectories of $f_\veps$ is reduced by the factor $r_0 r_1
\cdots r_{n-2}$.  Due to our definition of an admissible sequence $\r$,
the remaining factor $r_{n-1} r_{2n-1}$ does not affect the necessary
summability properties for the bounds
$\mu_n(C,\dt,\rho,\r,M_{1+\rho})$.  There is an additional
distortion of our estimates that is exponential in $n$, due to the fact that
an image of a finite-dimensional brick of $\eps$-parameters
under the Lagrange map is a parallelepiped of $u$-parameters, but
no longer a brick.  This exponential factor is also not problematic,
because our bound $\mu_n(C,\dt,\rho,\r,M_{1+\rho})$ decays
superexponentially in $n$.

\section{Conclusion}

In this announcement we have only been able to outline some of the
fundamental tools that are needed for the proof of the main result,
which will appear in \cite{K5} and \cite{KH}.  Here we list some of major
difficulties appearing in the proof.

$\bullet$ We must handle almost periodic trajectories of length $n$
that have a close return after $k<n$ iterates, so that as discussed
above the product of distances along the trajectory is small.  The
precise definition of a close return is a major problem here.  It must
not be too restrictive, because we must also show that a trajectory
without close returns is simple (the product of distances is not too
small).\footnote{This is exactly the place in the proof where we need to
impose superexponential decay of our bounds on hyperbolicity and
periodicity.}

$\bullet$ In dimension $N>1$, the Lagrange interpolation
polynomials involve products of differences of coordinates
of points, which may be small even though the points themselves
are not close.  Thus we must be careful about how we construct
the Lagrange basis for a given $n$-tuple of points
$x_0,\dots,x_{n-1} \subset B^N$ and how to incorporate this into
the general framework of the space of Lagrange interpolation polynomials.

$\bullet$ At $n$-th stage of the induction we need to deal with the
$(2n)^N$-dimensional space $W_{\leq 2n-1,N}$ of polynomials of degree
$2n-1$ in $N$ variables and handle the distortion properties of
the Lagrange map. In such a large dimensional space, even the ratio of
volumes of the unit ball and the unit cube is of order
$(2n)^{N(2n)^N}$ \cite{San}.

In \cite{K5}, \cite{KH}, based on \cite{K4}, we first prove the main
$1$-dimensional result for the case $N = \rho = 1$, discussed in
Sections~5, 6, and 7 of this announcement, and then using additional
tools and ideas complete the proof in the general case.

\section*{Appendix: \textup{Diff}$^r(B^N)$ and \textup{Diff}$^r(M)$}

Given a smooth ($C^\infty$) compact manifold $M$ of dimension $D$,
for $N > 2D$ the Whitney Embedding Theorem says that a generic smooth
function from $M$ to $\R^N$ is a a diffeomorphism between $M$ and its
image.  To simplify notation, we identify $M$ with its image, so that
$M$ becomes a submanifold of $\R^N$.

Let $U \subset \R^N$ be a closed neighborhood of $M$, chosen
sufficiently small that there is a well-defined projection $\pi : U
\to M$ for which $\pi(x)$ is the closest point in $M$ to $x$.  Then
for each $y \in M$, $\pi^{-1}(y)$ is an $(N - D)$-dimensional disk.
For $0 < \rho < 1$ and $y, z \in M$ choose a linear function
$g_{\rho,y,z} : \pi^{-1}(y) \to \pi^{-1}(z)$ that maps $y$ to $z$ and
contracts distances by a factor of $\rho$, and such that the
dependence of $g_{\rho,y,z}$ on $y$ and $z$ is $C^r$. Then we can
extend each $f \in \textup{Diff}^r(M)$ to a function $F \in C^r(U)$
that is a diffeomorphism from $U$ to a subset of its interior by
letting
$$
F(x) = g_{\rho,\pi(x),f(\pi(x))}(x)
$$
where $\rho = \|f^{-1}\|_{C^1}^{-r}/2$. Then by Fenichel's
Theorem \cite{F}, every sufficiently small
perturbation $F_\varepsilon \in C^r(U)$ of such an $F$ has an
invariant manifold $M_\varepsilon \subset U$ for which
$\pi|_{M_\varepsilon}$ is a $C^r$ diffeomorphism from $M_\varepsilon$
to $M$. Then to such an $F_\varepsilon$ we can associate
a diffeomorphism $f_\varepsilon \in \textup{Diff}^r(M)$ by letting
$$
f_\varepsilon(y) = \pi(F_\varepsilon(\pi|_{M_\varepsilon}^{-1}(y))).
$$

Notice that the periodic points of $F_\varepsilon$ all lie on
$M_\varepsilon$ and are in one-to-one correspondence with the periodic
points of $f_\varepsilon$.  Furthermore, because $f_\varepsilon$ and
$F_\varepsilon|_{M_\varepsilon}$ are conjugate, the hyperbolicity of
each periodic orbit is the same for either map.  Thus any estimate on
$P_n(F_\varepsilon)$ or $\gamma_n(F_\varepsilon)$ applies also to
$f_\varepsilon$.

\end{document}